\documentclass[11pt]{article}
\usepackage{authblk}

\usepackage[latin1]{inputenc}
\usepackage[T1]{fontenc}
\usepackage{amsmath}
\usepackage{amsfonts}
\usepackage{amssymb}

\usepackage[normalem]{ulem} 

\usepackage{amsmath,xcolor,ulem}
\usepackage{amsfonts}
\usepackage{amssymb,todonotes}
\usepackage{amsthm}
\usepackage{graphicx}
\usepackage{marginnote}

\usepackage[colorlinks=true, allcolors=blue]{hyperref}

\usepackage[top=2.8cm,bottom=2.8cm,left=2.5cm,right=2.5cm]{geometry}
\usepackage{float}
\usepackage[algoruled]{algorithm2e}
\usepackage[font=footnotesize,width=\linewidth]{caption}

\newtheorem{defi}{Definition}
\newtheorem{thm}{Theorem}
\newtheorem{prop}{Proposition}

\title{Integrating Port-Hamiltonian Systems with Neural Networks: \\ From Deterministic to Stochastic Frameworks}

\author{Luca Di Persio\footnote{ \href{mailto:luca.dipersio@univr.it}{ luca.dipersio@univr.it}} ,
Matthias Ehrhardt\footnote{Corresponding author, \href{mailto:ehrhardt@uni-wuppertal.de}{ehrhardt@uni-wuppertal.de}} ,
Sofia Rizzotto\footnote{ \href{mailto:sofia.rizzotto@studenti.univr.it}{sofia.rizzotto@studenti.univr.it}} 
}

\affil{Department of Computer Science -- College of  Mathematics\\ University of Verona, Italy}
\affil{IMACM, School of Mathematics and Natural Sciences, \\ University of Wuppertal, Germany}

\begin{document}
\maketitle

\begin{tikzpicture}[remember picture,overlay]
	\node[anchor=north east,inner sep=20pt] at (current page.north east)
	{\includegraphics[scale=0.2]{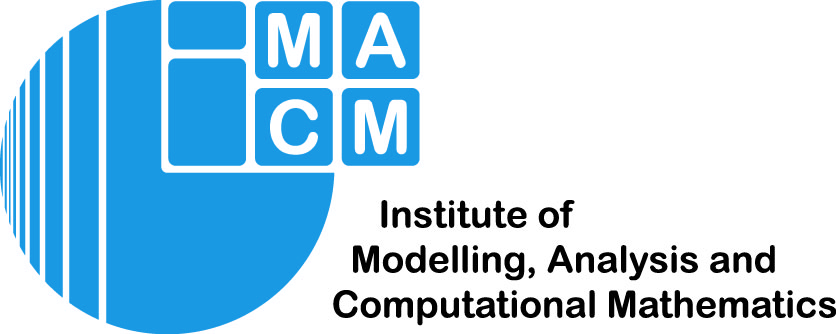}};
\end{tikzpicture}

\begin{abstract}
This article presents an innovative approach to integrating port-Hamiltonian systems with neural network architectures, transitioning from deterministic to stochastic models. 
The study presents novel mathematical formulations and computational models that extend the understanding of dynamical systems under uncertainty and complex interactions. 
It emphasizes the significant progress in learning and predicting the dynamics of non-autonomous systems using port-Hamiltonian neural networks (pHNNs). 
It also explores the implications of stochastic neural networks in various dynamical systems.
\end{abstract}

\begin{minipage}{0.9\linewidth}







\medskip

\end{minipage}

\section{Introduction}\label{sec1}
This work presents an innovative integration of port-Hamiltonian systems (PHS) with neural network architectures,
focusing in particular on the transition from deterministic models to those that incorporate stochastic elements. 
We first explore the advanced mathematical formulations and computational models to deepen the understanding of dynamical systems under uncertainty,
including complex interactions and measurement errors. 
We first recall the definition of port-Hamiltonian systems as a synthesis of port modeling and geometric Hamiltonian dynamics, 
highlighting the use of the Dirac structure, which generalizes the concept of Poisson and pre-symplectic structures, to represent the energetic topology of a system.
 We then consider the stochastic extension of port-Hamiltonian systems to analyze uncertainties such as measurement errors and unknown environmental effects, 
 making noise an intrinsic aspect of each port.
 For completeness, we also provide a detailed explanation of deterministic port-Hamiltonian systems using a coordinate-free geometric formulation. 
 
 The systems are described by generalized Poisson brackets and Hamiltonian functions, 
 which are translated into equations involving Jacobian matrices and vector fields associated with Hamiltonian functions in local coordinates.
 We then consider the associated extension to implicit systems, which includes systems with algebraic constraints. 
These take into account the state space representing stored energy, vector spaces of flow and effort variables representing power ports, 
and the Dirac structure outlining the power-saving connections within the system, emphasizing a significant property of these systems: passivity.
 Special emphasis is placed on its relationship to the properties of the Dirac structure and its importance in ensuring stability in control systems. 
 The transition from deterministic to stochastic \textit{port-Hamiltonian neural networks} (pHNNs) is derived by incorporating noise directly into the structure via ports. 
This initiative represents a critical conceptual innovation and is motivated by real applications where uncertainties and the need to accommodate errors affect systems.

A port-Hamiltonian system combines aspects of port modeling with a coordinate-free geometric Hamiltonian dynamics, 
which describes the dynamics of the physical system and the interconnection structure through the equations of motion. 
The main feature of these systems lies in the geometric structure defining them, the so-called Dirac structure,
which generalizes the concept of (pseudo) Poisson and pre-symplectic structures and represents the energetic topology of the system. 
This feature allows the framework to gain importance in modeling complex engineering systems with multiple energy exchange domains. 
The main focus of this paper is to explore an extension of port-Hamiltonian systems that incorporates random perturbations to account for uncertainties, 
measurement errors, and environmental interactions in dynamic systems. 
In particular, we will explore the innovation presented by Cordoni, Di Persio, and Muradore \cite{[3]:}, 
where the uniqueness lies in considering noise as an intrinsic aspect of each port. 
This is motivated by the need to account for measurement errors, parameter estimation uncertainties, and the unknown effects of the environment on the system.

The paper is organized a s follows.
Section~\ref{sec1}
establishes the basic concepts and explores different types of port-Hamiltonian systems, including discrete ones, and their application in neural networks.
Section~\ref{sec11} lays the foundation by introducing essential concepts such as explicit input-state-output port-Hamiltonian systems, local coordinates, and implicit port-Hamiltonian systems. 
Next, Section~\ref{sec12} delves into a special type of port-Hamiltonian systems: discrete ones. 
It covers the generalized Dirac structure, a mathematical construct that underlies these systems, and how two such structures can be connected. 
%
Section~\ref{sec2} 
focuses on extending the concept of port-Hamiltonian systems to include random elements, making them suitable for modeling systems with uncertainties or randomness. We generalize the concept of Dirac structure to stochastic ones, 
to account for the stochastic nature introduced in this section.
    We also introduce the concept of a noise port, a special type of port through which random fluctuations or uncertainties enter the system.
%
    Section~\ref{sec3} 
    extends the concept of passivity developed for deterministic systems to include the presence of randomness. Hereby, we explain strong and weak passivity and show 
    how to determine whether a given stochastic port-Hamiltonian systems (SPHS) exhibits passivity.
    %
    Section~\ref{sec4} 
    focuses on applying and extending the concepts of stochastic PHS to various scenarios. 
    First, Section~\ref{sec41} explores how to connect and analyze multiple SPHS systems.
    Section~\ref{sec42} 
    adapts the continuous-time framework of SPHS to represent systems that evolve in discrete time steps.
 Section~\ref{sec43} 
discusses the application of SPHS to model the stochastic motion of agents, e.g.\ cars. 

In Section~\ref{sec5} 
show successful applications of pHNNs in tasks such as simulating a damped mass-spring system and a chaotic Duffing system. 
These examples show that pHNNs can not only learn system dynamics, but also recover complex behaviors such as chaotic trajectories with minimal data.
Overall, pHNNs hold promise for several areas involving complex physical systems, including chemical interactions, robot motion control, and understanding general system dynamics without requiring precise details of the forces involved.

In the conclusions in Section~\ref{sec6}, we summarize our work, which is based on the method of Colonius and Gr\"une \cite{[53]:}, which used neural networks for controller design. 
Our focus is on adapting this method to handle stochastic port-Hamiltonian systems with random fluctuations with the two key steps:
1) Modeling the system with noise: Introducing stochastic differential equations (SDEs) that include a term representing random noise ($dW$) to describe the system dynamics. \linebreak
2)~Design a Robust Controller: Modify the control law to account for the noise. This involves adding a new term ($\mu(x)$) to the original control law ($K(x)$) to counteract the stochastic disturbances and ensure that the system remains stable and performs well.

Finally, in \ref{appA} we provide a discussion of stochastic neural networks (SNNs), which are inspired by the biological brain. 
Unlike traditional deterministic networks, SNNs incorporate randomness to improve training and avoid overfitting.




\subsection{Basic Concepts}\label{sec11}
The description given by Cordoni, Di Persio and Muradore \cite{[3]:} of an input-state-output deterministic port-Hamiltonian system
can be done by considering a geometric coordinate-free formulation in terms of Poisson brackets, i.e.\ the following system
\begin{equation}\label{a}
    \begin{cases}
        \dot{x} &= X_H(x)+\sum^m_{i=1} u_i X_{H_{g_i}}(x),\\
            y_i &= \{H,H_{g_i}\},
    \end{cases}
\end{equation}
is called a \textit{(explicit) input-state-output port-Hamiltonian system} (PHS) on a Poisson manifold $(X , \{\cdot,\cdot\})$ 
with a Hamiltonian function $H\in C^\infty(\mathcal{X})$, $x\in\mathbb{R}^n$, the $i$-th input $u_i\in U$, the $i$-th output $y_i\in U^*$ 
and the Hamiltonian vector field $X_{H_{g_i}}$ associated with the Hamiltonian $H_{g_i}$. In local coordinates, the previous system reads
\begin{equation}
    \begin{cases}
        \dot{x} &= J(x)\partial_xH+\sum_{i=1}^m u_i g_i(x),\\
           y_i  &= g_i^\top(x) \partial_xH.
    \end{cases}
\end{equation}
Moreover, considering $X_H^L(\cdot):=[\cdot,H]_L$, we can define the \textit{(explicit) input-state-output port-Hamiltonian system with dissipation} as
\begin{equation}\label{b}
    \begin{cases}
        \dot{x}&=X_H^L(x)+\sum_{i=1}^mu_iH_{g_i}(x),\\
        y_i&=[H,H_{g_i}],
    \end{cases}
\end{equation}
and in local coordinates this reads
\begin{equation}
    \begin{cases}
        \dot{x} &= \bigl(J(x)-R(x)\bigr) \partial_xH(x) + \sum_{i=1}^m u_i g_i(x),\\
            y_i &= g_i^\top(x) \partial_xH(x),
    \end{cases}
\end{equation}
where $R(x):=(g^R(x))^\top \Tilde{R}(x) g^R(x)$. 

Let's consider a physical system consisting of elements that store energy, a set of elements that dissipate energy, 
and a set of power ports interconnected by power-preserving links that can only transfer energy, not produce it. 
We can describe such a system by extending the framework of port-Hamiltonian systems to the context of implicit systems, 
i.e., systems with algebraic constraints.
Given a state space $\mathcal{X}$ (a smooth manifold whose elements represent the energy stored in the system), 
a vector space of flow variables $\mathcal{V}$ and its dual space of effort variables $\mathcal{V}^*$ (representing the power ports), 
a geometric Dirac structure $\mathcal{D}$ and a Hamiltonian function $\mathcal{H}$ representing the total energy of the system in a given state,
we can define an \textit{implicit port-Hamiltonian system} corresponding to $(\mathcal{X},\mathcal{V},\mathcal{D},\mathcal{H})$ as
\begin{equation}
    v=-\dot{x}\quad\text{and}\quad v^*=\frac{\partial\mathcal{H}}{\partial x}(x),
\end{equation}
which means that the system is defined by
\begin{equation}
    \Bigl(-\dot{x},\frac{\partial H}{\partial x}(x),f,e\Bigl)\in\mathcal{D}(x).
\end{equation}
Note that the above Dirac structure outlines the behavior of the internal connection, 
and its main property is that the power-conserving combination of Dirac structures remains a Dirac structure. 
This implies that any power-conserving connection of port-Hamiltonian systems is also a port-Hamiltonian system, 
where the Dirac structure is the composition of the Dirac structures of its elementary parts, and the Hamiltonian is the sum of the Hamiltonians.

Another important property of a system is passivity, which means that the energy input to the system is always greater than or equal to (lossless case) 
the energy output from the system. It is crucial for ensuring stability in control systems. 
Passivity is a consequence of the properties of the Dirac structure and the energy-dissipation relation for port-Hamiltonian systems.

\subsection{Discrete Systems}\label{sec12}
\begin{figure}[htbp]
\centering
\setlength{\unitlength}{0.8cm}
\includegraphics[width=.90\textwidth]{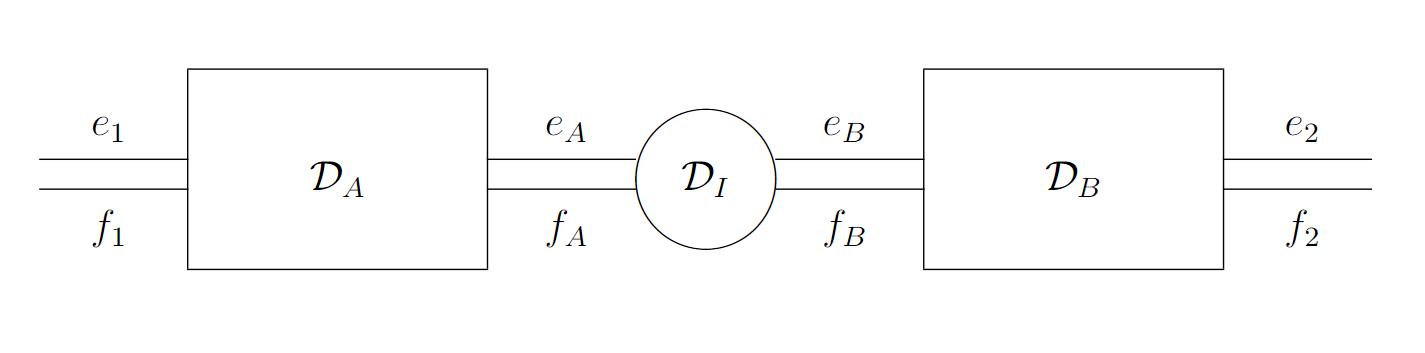}
\caption{Interconnection of two port-Hamiltonian systems.}\label{figure12}
\end{figure}

So far, we have described port-Hamiltonian systems in continuous time. 
However, it is essential to investigate whether these properties persist when we discretize time. 
As suggested by Viswanath, Clemente-Gallardo and van der Schaft \cite{[27]:}, the discretization of a Hamiltonian system can be done while preserving energy conservation. 
Discrete systems can be obtained by discretizing a continuous system or by modelling systems directly at the discrete level. 
To obtain the latter, we can use Poisson brackets, since they satisfy the properties of skew symmetry, bilinearity, and a modified Leibniz rule, 
which are necessary to preserve the structure of the Hamiltonian dynamics in the discrete setting. 
First, we need to define the discrete Dirac structure. 
For this purpose we denote with $\mathfrak{X}(A)$ a space of discrete vector fields, with $\Lambda^1(A)$ the space of discrete 1-forms. 
Then a \textit{generalized Dirac structure} on an $n$-dimensional discrete manifold 
is an $n$-dimensional linear subspace $\mathcal{D}\subset\mathfrak{X}(A)\times\Lambda^1(A)$ such that $\mathcal{D}=\mathcal{D}^\top$ with
\begin{equation*}
    \mathcal{D}^\top=\bigl\{(Y,\beta)\in\mathfrak{X}(A)\times\Lambda^1(A)\quad
    \text{where}\quad\langle \alpha,X\rangle+\langle\beta,Y\rangle=0,\quad
    \forall\,(X,\alpha)\in\mathcal{D}\bigr\},
\end{equation*}
where $\langle\cdot,\cdot\rangle$ is the pairing between $\mathbb{F}^n$ and $\mathbb{F}^{n*}$. 
Then, assume to have the situation illustrated in Figure~\ref{figure12}, 
denoting with $\mathcal{F}_i$ the space of flow and with $\mathcal{F}_i^*$ the space of effort of the Dirac structure $\mathcal{D}_i$ with $i=A,B$, 
we can define the \textit{interconnection between the two Dirac structures} $\mathcal{D}_A$ and $\mathcal{D}_B$ as
\begin{equation}
    \begin{split}
        \mathcal{D}_A\circ\mathcal{D}_B &:= \big\{(f_1,e_1,f_2,e_2)\in\mathcal{F}_1\times\mathcal{F}_1^*\times\mathcal{F}_2\times\mathcal{F}_2^*\quad\text{such that}\\ 
        &\exists(f,e)\in\mathcal{F}\times\mathcal{F}^*\quad\text{with}\quad
        (f_1,e_1,f,e)\in\mathcal{D}_A\quad\text{and}\quad(-f,e,f_2,e_2)\in\mathcal{D}_B\big\}.
        \end{split}
\end{equation}
Then $\mathcal{D}_A\circ\mathcal{D}_B$ is a Dirac structure.

\begin{defi}[Implicit discrete port-Hamiltonian system \cite{[27]:}]
    Let $H:\mathcal{Z}\to\mathbb{F}$ be a discrete Hamiltonian, $\mathcal{F}_P$ be the space of external flows $f$, $\mathcal{E}_P=\mathcal{F}_P^*$ the space of external effort $e$, $\mathcal{D}$ be the Dirac structure which depends only on the coordinate $z$, then the \textit{implicit discrete port-Hamiltonian system} is defined as
    \begin{equation}
        \Bigl(-\frac{\Delta z}{\Delta t},f,\Game_zH(z),e\Bigl)\in\mathcal{D}(z).
    \end{equation}
\end{defi}

There are also some cases where the system consists of an interconnection of smooth and discrete port-Hamiltonian systems, 
which can be achieved (as proposed by Viswanath, Clemente-Gallardo and van der Schaft \cite{[30]:})
by transforming the problem of energy-conserving connection on sampling instances, into an energy-conserving interconnection with an external flow source. 
In this way, we ensure energy conservation over sampling intervals. 
Moreover, the resulting interconnected system is again passive and port-Hamiltonian. 
Thus, the integration of discrete computational models with physical systems does not compromise the stability and performance characteristics of passive systems. 
However, energy conservation does not hold, so the idea of Kotyczka and Lefevre \cite{[28]:} is to introduce a discrete-time Dirac structure and a discrete-time port-Hamiltonian system.
This method approximates the continuous-time structural energy balance and uses symplectic numerical time integration by collocation methods.


\section{Stochastic port-Hamiltonian systems}\label{sec2}
Consider a complete probability space $\bigl(\Omega,\mathcal{F},(\mathcal{F}_t)_{t\in\mathbb{R}_+},\mathbb{P}\bigl)$ 
and denote by $\delta Z$ the Stratonovich integration and by $dZ$ the It\^o integration along the semimartingale $Z$. 
As mentioned above, the stochastic port-Hamiltonian formulation allows for a more extensive source of randomness by treating each element of the system as a semimartingale. 
The adoption of the Stratonovich calculus is due to the geometric nature of the Dirac structure, 
but we can translate it into the corresponding It\^o formulation to exploit its probabilistic properties. 
Consider the following system
\begin{equation}\label{lcsphs}
    \begin{cases}
        \delta X_t &= \bigl(J(X_t)-R(X_t) \bigr)\partial_xH(X_t)\delta Z_t 
        + g(X_t)u\delta Z^g_t + \xi(X_t)\delta Z^N_t,\\
               y_t &= g^\top(X_t)\partial_xH(X_t),
    \end{cases}
\end{equation}
where $R(x):=(g^R(x))^\top\Tilde{R}(x)g^R(x)$, $W$ is a Brownian motion, and $Z$, $Z^g$, and $Z^N$ are semimartingales. 
Then \eqref{lcsphs} describes the stochastic port-Hamiltonian system in local coordinates.

By generalizing the notion of Dirac structure, it is possible to extend the previously 
described framework by incorporating scenarios where noise is introduced into the system as a stochastic external random field
and as a random perturbation of any port connected to the system.
\begin{defi}[Orthogonal complement \cite{[3]:}]
    Given a manifold $\mathcal{X}$, $I\subset\mathbb{R}_+$, a bundle $\mathcal{D}\subset T\mathcal{X}\oplus T^*\mathcal{X}$, a differential 1-form $\sigma$ on $\mathcal{X}$ and an integral curve $X:I\to\mathcal{X}$ of a Stratonovich vector field $\delta X_t$, the \textit{orthogonal complement} of $\mathcal{D}$ is 
    \begin{equation}
        \mathcal{D}^\perp=\bigl\{(\delta X_t,\sigma)\subset T\mathcal{X}\oplus T^*\mathcal{X}:\int_0^t\langle \sigma,\delta \bar{X}_s\rangle + \int_0^t\langle \bar{\sigma},\delta X_s\rangle =0,\; 
        \forall\,(\delta\bar{X}_t,\bar{\sigma})\in\mathcal{D},\ t\in I\bigr\}
    \end{equation}
\end{defi}
\begin{defi}[Generalized Dirac structure \cite{[3]:}]
With the same notation as above, we call \textit{generalized stochastic Dirac
structure} a smooth vector subbundle $\mathcal{D}\subset T\mathcal{X}\oplus T^*\mathcal{X}$ such that $\mathcal{D}=\mathcal{D}^\perp$.
\end{defi}
\begin{defi}[Implicit generalized stochastic PHS \cite{[3]:}]
    Let $H:\mathcal{X}\to\mathbb{R}$ be a Hamiltonian function, 
    $Z$ a semimartingale perturbing the system, then an \textit{implicit generalized stochastic port-Hamiltonian system}
    on $\mathcal{X}$ is a 4-tuple $(\mathcal{X},Z,\mathcal{D},H)$ such that
    \begin{equation}
        \bigl(\delta X_t,\textbf{d}H(X_t)\bigr)\in\mathcal{D}(X_t)\quad 
        \forall\,t\in I.
    \end{equation}
        Including a resistive element and an external element control then an \textit{implicit generalized port-Hamiltonian system with resistive structure} $\mathcal{R}$ is a 5-tuple $(\mathcal{X},\textbf{Z},\mathcal{F},\mathcal{D},H)$ such that
    \begin{equation*}
        (-\delta X_t,\textbf{d}H,\delta f_t^R,e_t^R,\delta f_t^C,e_t^C)\in\mathcal{D(X)_t}\quad
        \text{with}\quad (\delta f_t^R,e_t^R)\in\mathcal{R(X)}_t.
    \end{equation*}
\end{defi}
\begin{prop} \cite{[3]:}
    Implicit port-Hamiltonian systems satisfy energy conservation property that is
    \begin{equation}\label{cons}
        H(X_t)-H(X_0)=\int_0^t\langle \textbf{d}H,\delta X_s\rangle 
    \end{equation}
        which can be written in short notation as
    \begin{equation}
        \delta H(X_t)=\langle \textbf{d}H,\delta X_t\rangle .
    \end{equation}
    The energy balance is
    \begin{equation}
        H(X_t)-H(X_0)=\int_0^t\langle e_s^R,\delta f_s^R\rangle +\int_0^t\langle e_s^C,\delta f_s^C\rangle \le\int_0^t\langle e_s^C,\delta f_s^C\rangle .
    \end{equation}
However, the condition $\int_0^t\langle e_s^R,\delta f_s^R\rangle \le0$ imposed on the resistive port is too strong since it is difficult for
it to happen in practice, so the idea is to introduce the following weaker definition of the resistive relation $\mathcal{R}_W\subset\mathcal{F}_{Z_R}\times\mathcal{E}_R$:
\begin{equation}
    \mathbb{E}\int_0^t\langle e_s^R,\delta f_s^R\rangle \le0.
\end{equation}
Consequently, the \textit{mean power balance} requires that the energy be conserved and dissipated in mean value, i.e.\ it reads
\begin{equation}
    \mathbb{E}\bigl(H(X_t)-H(X_0)\bigr) = \mathbb{E}\int_0^t\langle e_s^R,\delta f_s^R\rangle +\mathbb{E}\int_0^t\langle e_s^C,\delta f_s^C\rangle \le\mathbb{E}\int_0^t\langle e_s^C,\delta f_s^C\rangle .
\end{equation}
\end{prop}

\begin{figure}[htbp]
\centering
\setlength{\unitlength}{0.6cm}
\includegraphics[width=.90\textwidth]{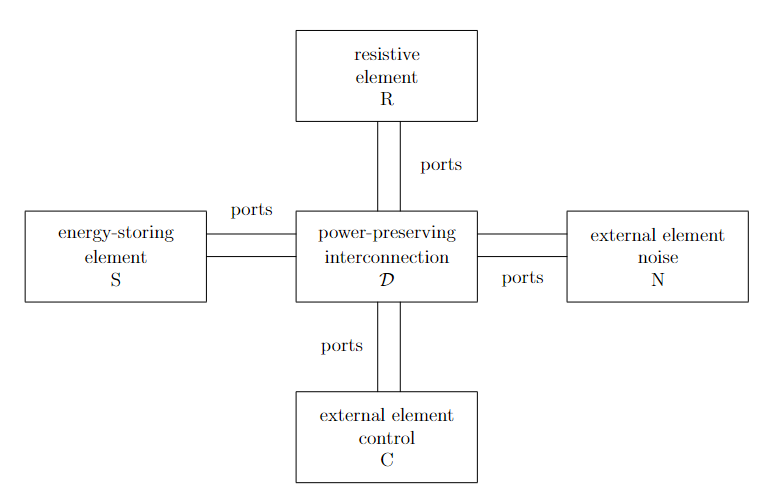}
\caption{Schematic representation of a general implicit port-Hamiltonian system.}\label{figure10}
\end{figure}

It is possible to further generalize the Hamiltonian by introducing an external perturbation of the system,
i.e.\ a new type of port called \textit{noise port} perturbed by the semimartingale $Z_N$ (see Figure~\ref{figure10}). 
Thus, in this case the \textit{implicit generalized stochastic port-Hamiltonian system with resistive structure} is a 5-tuple $(\mathcal{X},\textbf{Z},\mathcal{F},\mathcal{D},H)$ such that
    \begin{equation}\label{64}
        \bigl(-\delta X_t,\textbf{d}H,\delta f_t^R,e_t^R,\delta f_t^C,e_t^C,\delta f_t^N,e_t^N\bigr)\in\mathcal{D}(X_t)
    \end{equation}
    and the weak energy balance is given by 
    \begin{equation}
        \mathbb{E}H(X_t)-\mathbb{E}H(X_0)\le\mathbb{E}\int_0^t\langle e_s^N,\delta f_s^N\rangle +\mathbb{E}\int_0^t\langle e_s^C,\delta f_s^C\rangle.
    \end{equation}

Furthermore, we can reformulate the stochastic port-Hamiltonian system in It\^o form by defining 
the following alternative representation of the Dirac structure.
\begin{defi}[Dirac structure \cite{[3]:}]
    Let $\mathcal{F}:=\mathcal{F}_{Z_R}\times\mathcal{F}_{Z_C}\times\mathcal{F}_{Z_N}$ be the space of
    flows $\delta f$, $\mathcal{E}=\mathcal{F}^*$ be the dual space of efforts $e$, $G_\theta:\mathcal{F}_{Z_\theta}\to\mathcal{F}_{Z_\theta}$ be a function such that $\langle e_t^S,G_\theta\delta f_t^\theta\rangle =\langle G_\theta^*e_t^S,\delta f_t^\theta\rangle$ and $J$ be a matrix such that $J=-J^\top$. 
    Then the Dirac structure $\mathcal{D}$ can be defined as
    \begin{equation}
    \begin{split}
        \mathcal{D}:=\bigl\{(&\delta f_t^S,\delta f_t^R,\delta f_t^C,\delta f_t^N,e_t^S,e_t^R,e_t^C,e_t^N)\in\mathcal{F}\times\mathcal{E}:\\
        &\delta f_t^S=-Je_t^S\delta Z_t-G_R\delta f_t^R-G_C\delta f_t^C-G_N\delta f_t^N,\\
        &e_t^R=G_R^*e_t^S,\quad 
        e_t^C=G^*_Ce_t^S,\quad 
        e_t^N=G^*e_t^S\bigr\}.
    \end{split}
    \end{equation}
\end{defi}
Consider the special case 
\begin{equation}
\begin{split}
    &\delta f_t^R=-\Tilde{R}e_t^R\delta Z_t,\quad 
    \delta f_t^N=\xi_t\delta Z_t^N,\quad 
    \delta f_t^C=u_t\delta Z_t^C\\ &\text{with} \quad 
    \mathbb{E}\int_0^t\langle e_s^R,\Tilde{R}e_s^R\delta Z_s\rangle -\mathbb{E}\int_0^t\langle e_s^N,f_s^N\delta Z_s^N\rangle \ge0,
\end{split}
\end{equation}
where the reason for the minus sign in front of $\Tilde{R}$ is that we want it to be the incoming power regarding the interconnection (as in \cite{[2]:}).
\begin{defi}[Stochastic input-output PHS with stochastic Dirac structure \cite{[3]:}]
Using the same notation as above, if $\textbf{Z}=(Z,Z^R,Z^C,Z^N)$ is a semimartingale and $H:\mathcal{X}\to\mathbb{R}$ is a Hamiltonian function, then the \textit{stochastic input-output port Hamiltonian system with stochastic Dirac structure} is given by
    \begin{equation}\label{mio}
        \begin{cases}
            \delta X_t &= -JdH(X_t)\delta Z_t+G_R\Tilde{R}e_t^R\delta Z_t-G_Cu_t\delta Z_t^C-G_N\xi_t\delta Z_t^N,\\
            e_t^N &= G_N^*\textbf{d}H(X_t),\\
            e_t^C &= G_C^*\textbf{d}H(X_t)
        \end{cases}
    \end{equation}
    and by taking $\Tilde{J}=-J$ and $e_t^R=G^*_R\mathbf{d}H(X_t)$, the system \ref{mio} becomes
     \begin{equation}\label{nonmio}
        \begin{cases}
            \delta X_t &= \big(\Tilde{J}+G_R\Tilde{R}G^*_R\big)\textbf{d}H(X_t)\delta Z_t-G_Cu_t\delta Z_t^C-G_N\xi_t\delta Z_t^N,\\
            e_t^N &= G_N^*\textbf{d}H(X_t),\\
            e_t^C &= G_C^*\textbf{d}H(X_t).
        \end{cases}
    \end{equation}
\end{defi}
\begin{thm} \cite{[3]:}\label{conv}
    If $X$ is a solution of the equation \eqref{nonmio} and $Z,Z^N,Z^C$ are such that 
    \begin{equation*}
        \langle Z,Z^C\rangle _t = \langle Z,Z^N\rangle _t = \langle Z,Z^C\rangle _t=0
    \end{equation*}
    where $\langle \cdot,\cdot\rangle_t$ is the quadratic covariation at time $t$, 
    then $X$ can be equivalently rewritten in It\^o terms as
    \begin{equation}
        \begin{split}
            dX_t &= V^S(X_t)\,dZ_t + \mathcal{L}_{V^S}V^S(X_t)\,d\langle Z,Z\rangle_t+\\
            &\qquad -\sum_{i=1}^{n^N} V_i^N(X_t)\,dZ_t^N - \frac{1}{2}\sum_{i,j=1}^{n^N}\mathcal{L}_{C_j^N} V_i^N(X_t)\,d\langle Z^{N;i},Z^{N;j}_t\rangle +\\
            &\qquad -\sum_{i=1}^{n^C}V_i^C(X_t) u_t^i\,dZ_t^{C;i} -\frac{1}{2}\sum_{i,j=1}^{n^C}\mathcal{L}_{V_j^C}V_i^C(X_t) u_t \,d\langle Z^{C;i},Z^{C;j}\rangle_t,
        \end{split}
    \end{equation}
    where $\mathcal{L}$ is the Lie derivative and $V^\alpha$, $\alpha=S,N,C$ are defined as
    \begin{equation}\label{V}
    \big(\Tilde{J}+G_R\Tilde{R}G^*_R\big)\,\textbf{d}H=V^S,\quad
    G_N\xi_t=\sum_{i=1}^{n^N}V_i^M,\quad
    G_C=\sum_{i=1}^{n^C}V^C_i.
\end{equation}
\end{thm}

\section{Passivity in Stochastic Systems}\label{sec3}
Extending the notion of passivity to the stochastic case is challenging because the presence of noise affects the energy of the system. 
The standard requirement that the structure matrix $R$ be symmetric and positive semidefinite is no longer sufficient 
due to the presence of the semimartingale $Z$, which makes the Stochastic port-Hamiltonian System (SPHS) non-dissipative,
 so it is necessary to impose specific conditions on the noise to ensure losslessness and passivity.
\begin{prop} \cite{[3]:}
    If $X$ is the solution of an explicit I-S-O stochastic PHS with dissipation, i.e.
    \begin{equation}\label{eisosd}
        \begin{cases}
        \delta X_t &= X^L_H(X_t)\,\delta Z_t + uX^L_{H_g}(X_t)\,\delta Z^g_t + X^L_{H_N}(X_t)\,\delta Z^N_t,\\
        y_t &= [H,H_g]_L,
    \end{cases}
    \end{equation}
    with 
    \begin{equation*}
      X^L_H(\cdot) := [\cdot, H]_L,\quad 
      X^L_{H_g}(\cdot) := [\cdot, H_g]_L,\quad
      X^L_{H_N}(\cdot) := [\cdot, H_N ]_L,
    \end{equation*}
    then for all $\varphi\in C^\infty(\mathcal{X})$ it holds
    \begin{equation}\label{35}
        \begin{cases}
            \delta\varphi(X_t)&=[\varphi,H]_L(X_t)\,\delta Z_t + u[\varphi,H_g]_L(X_t)\,\delta Z_t^g+[\varphi,H_N]_L(X_t)\,\delta Z_t^N,\\
            y_t&=[H,H_g]_L.
        \end{cases}
    \end{equation}
\end{prop}
The equation changes when there is no external noise and the semimartingale perturbing the control is deterministic.
\begin{defi}[Strong and weak passivity \cite{[3]:}]
    If $H\in C^\infty(\mathcal{X})$ is the total energy of the explicit I-S-O stochastic PHS with dissipation,
    then it is \textit{strongly passive} if for all $t\ge0$ it holds
    \begin{equation}\label{sp}
        H(X_t)\le H(X_0) + \int_0^t u^\top(s) y(s) \,\delta Z_s^C,
    \end{equation}
    or \textit{weakly passive} if for all $t\ge0$ it holds
    \begin{equation}\label{wp}
        \mathbb{E}H(X_t)\le\mathbb{E}H(X_0) + \mathbb{E}\int_0^t u^\top(s) y(s) \delta Z_s^C.
    \end{equation}
\end{defi}
Assuming $\xi=0$, the energy conservation relation of the system 
\begin{equation}\label{81}
        \begin{cases}
            \delta X_t &= \big(J(X_t)-R(X_t)\big)\partial_xH(X_t)\,\delta Z_t + g(X_t)u\,\delta Z_t^C + \xi(X_t)\,\delta Z_t^N,\\
            e^C &= g^\top(X_t) \partial_x^H(X_t),
        \end{cases}
    \end{equation}
where $J=J^\top$, is given by
\begin{equation}\label{calc}
    \begin{split}
        H(X_t)-H(X_0) &= \int_0^t\langle \textbf{d}H,\delta X_s\rangle
        = \int_0^t \langle e_s^R,\delta f_s^R\rangle 
           +\int_0^t \langle e_s^C,\delta f_s^C\rangle \\
        &= \int_0^t \bigl\langle \partial_xH(X_s),R(X_s)e^R\,\delta Z_s \bigr\rangle +\int_0^t \langle y,u\,\delta Z_s^C\rangle \\
        &= \int_0^t \bigl\langle \partial_xH(X_s),-R(X_s)\partial_xH(X_s)\delta Z_s\bigr\rangle 
        + \int_0^t y^\top u\,\delta Z_s^C\\
        &=-\int_0^t \partial_x^\top H(X_s)R(X_s)\partial_xH(X_s)\delta Z_s+\int_0^ty^\top u\,\delta Z_s^C.
    \end{split}
\end{equation}
We can see that even if $R$ is strictly positive, 
we cannot infer that the strongly passive condition holds;
 in fact, we should also require that
\begin{equation}\label{adcond}
    \int_0^t \partial_x^\top H(X_s) R(X_s) \partial_xH(X_s) \,\delta Z_s\ge0.
\end{equation}
However, the latter condition is usually difficult to satisfy in a real system, so it makes more sense to require the weaker condition.
\begin{equation}
    \mathbb{E}\int_0^t\partial_x^\top H(X_s)R(X_s)\partial_xH(X_s)\,\delta Z_s\ge0,
\end{equation}
which guarantee the weak passivity property. However, computing the expectation of a Stratonovich integral is often challenging. 
The natural solution to this problem is to transform the equation \eqref{wp} into It\^o terms. 
This transformation allows us to take advantage of the favourable probabilistic properties of the It\^o integral.
 To achieve this, we can apply the Theorem~\ref{conv} to the equation \eqref{wp}.

\section{Application and Generalization}\label{sec4}
\subsection{Interconnection of Multiple SPHS}\label{sec41}
An essential property of port-Hamiltonian systems is their interconnectivity, which allows complex systems to be viewed as compositions of simpler parts. 
This interconnectivity can be analysed regarding the components and how they are interconnected.
 In particular, through the composition of Dirac structures, the power-preserving interconnection of port-Hamiltonian systems defines another port-Hamiltonian system. 
The Hamiltonian of the interconnected port-Hamiltonian system is the sum of the Hamiltonians of its components,
 and the energy-dissipation relation is the union of the energy-dissipation relations of the subsystems.

The following proposition discusses the connection of multiple SPHS, defining a new system with interconnected Dirac structures and combined Hamiltonians.

\begin{figure}[htbp]
\centering
\setlength{\unitlength}{0.8cm}
\includegraphics[width=.90\textwidth]{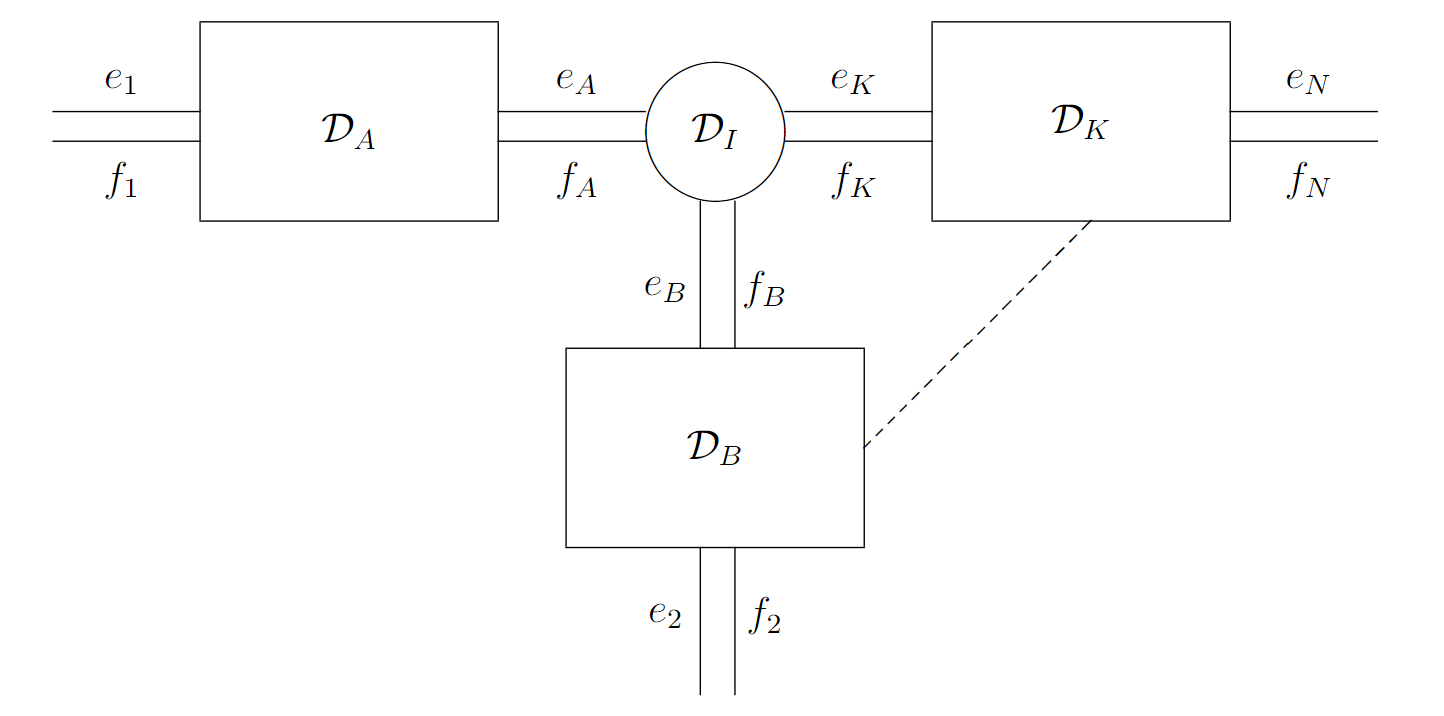}
\caption{Interconnection of $N$ implicit port-Hamiltonian systems.}\label{figure13}
\end{figure}

\begin{prop}
Suppose we have $N$ stochastic port-Hamiltonian systems with state space $\mathcal{X}_i$, Hamiltonian $H_i$,
flow-effort space $\mathcal{F}_i\times\mathcal{E}_i$ and perturbation $\textbf{Z}^i$ for $i=1,\dots,N$. 
Assuming that they are connected by $\mathcal{D}_I$ (see Figure~\ref{figure13}), 
then their interconnection defines a stochastic port-Hamiltonian system with Dirac structure $\mathcal{D}\circ\mathcal{D}_I$ and Hamiltonian $H:=\sum_{i=1}^NH_i$.
\end{prop}

\subsection{Discrete Stochastic PHS}\label{sec42}
Consider the continuous stochastic port-Hamiltonian system 
\begin{equation}\label{new}
    \begin{cases}
        dX_t &= \big((J-R)\partial_xH(X_t) + g(X_t)u_t\big)\,dt + \xi(X_t) \,\delta W_t,\\
        y_t &= g^\top(X_t)\partial_xH(X_t),\\
        z_t &= \xi^\top(X_t)\partial_xH(X_t),
    \end{cases}
\end{equation}
where $J =-J^\top$, R positive semidefinite, $g$ represents control port, $H$ is the Hamiltonian, $u\in U$ is the control
input, $y\in Y$ is the output of the system, $\xi$ is a matrix, $z$ is the associated effort to $\delta W_t$ and $W$ is a standard Brownian motion adapted to the reference filtration $(\mathcal{F}_t)_{t\ge0}$. Then we can introduce the discretization
\begin{equation}\label{def}
    \dot{X}(t_0^k + \tau h) = -f(t_0^k+\tau h) = -\sum_{j=1}^s f_j^k l_j(\tau)
\end{equation}
with
\begin{equation*}
    \dot{X}(t_i^k):=-f_i^k,\quad 
    l_i=\prod_{j=1}^s\frac{\tau-c_j}{c_i-c_j},\quad 
    \tau\in[0,1],
\end{equation*}
where $l_i$ is the $i^{\rm th}$ Lagrange interpolation polynomial of order $s$ and $\tau$ is the normalized time parameterizing the sampling intervals. 
Thus, we can generalize the continuous SPHS to a discrete form, preserving the structure of the Hamiltonian and the control inputs, as follows:
\begin{defi}[Discrete stochastic port-Hamiltonian system \cite{[28]:}]
    A \textit{discrete stochastic port-Hamiltonian system} can be written as
    \begin{equation}\label{dsphs}
        \begin{cases}
            X(t_0^k+c_ih) &= x_0^k - h\sum_{i=1}^sa_{ij}f_j^k,\\
            X(t_0^k+h) &= x_0^k - h\sum_{i=1}^sa_jf_i^k,\\
            -a_{ij}f^k &= (J_j^k-R_j^k)a_{ij}e^k_j+a_{ij}g_j^ku_j^k+b_{ij}\xi^k_j\,\Delta W,
        \end{cases}
    \end{equation}
    where $\Delta W$ is a truncated centred Gaussian random variable with variance $h$,
    $a_{ij}=\int_0^{c_i}l_j(\sigma)d\sigma$, 
    $a_j=\int_0^1l_j(\sigma)d\sigma$ and $M=M^\top$ \cite{[28]:}.
\end{defi}
Note that in the discrete case, the system is \textit{passive} if it holds
    \begin{equation}
        \mathbb{E}\bigl[\Delta H^k\bigr]
        \le h\,\mathbb{E}\bigl[(y^k)^\top u^k\bigr].
    \end{equation}

\subsection{Stochastic Motion Model of Agents}\label{sec43}
An application to a stochastic motion model of agents is presented by Ehrhardt, Kruse, and Tordeux \cite{[48]:}, 
who analyze the case where positions and velocities of agents are modeled in a ring structure. 
The initial positions and velocities are set, and the system is governed by differential equations involving the velocities and positions of neighboring agents \cite{[48]:}.
The dynamic equations of the agents are given by:
     \begin{equation}\label{motion}
    \begin{cases}
        dQ_n(t) &= \bigl(p_{n+1}(t)-p_n(t)\bigr)\,dt,\\
        dp_n(t) &= \bigl(U'(Q_n(t))-U'(Q_{n-q}(t))\bigr)\,dt+ \\ 
             &\qquad +\beta\bigl(p_{n+1}(t)-2p_n(t)+p_{n-1}(t)\bigr)\,dt + \sigma\,dW_n(t)
    \end{cases}
\end{equation}
with $Q(0)=Q_0\in[0,+\infty)^N$ the initial distance, $p(0)=p_0$ the initial velocity, $\beta\in(0,+\infty)$ a dissipation rate, $\sigma\in\mathbb{R}$ the noise volatility, 
$U'$ the derivative of a convex potential $U\in C^1(\mathbb{R},[0,+\infty))$ 
and $W=(W_n)_{n=1}^N\colon[0,+\infty)\times\Omega\to\mathbb{R}^N$ an $N$-dimensional standard Brownian motion defined on $(\Omega,\mathcal{F},\mathbb{P})$.
These equations represent the acceleration of the $n^{\rm th}$ agent depending on the velocities of its neighbors and the stochastic perturbations of the Brownian motion \cite{[49]:}.

The motion of the agents is further formulated using a stochastic port-Hamiltonian framework
\begin{equation*}
    dZ(t) = (J - R)\nabla H(Z(t))\,dt + G\,dW(t),
\end{equation*}
where $Z(t)=\big(Q(t),p(t)\big)^\top\in\mathbb{R}^{2N}$, $t\in[0,+\infty)$, $J$ and $R$ are defined as skew-sym\-me\-tric and symmetric positive semidefinite matrices, respectively. 
This formulation allows the application of Hamiltonian dynamics to model agents' behaviour under stochastic influences. 
In particular, the Hamiltonian is independent of $Q$ and $U$, and its expectation could increase with time. 
Moreover, describing the limiting behavior of these stochastic systems is challenging, so the authors \cite{[48]:} 
focus on the specific scenario where the quadratic function characterizes the potential
\begin{equation}
    U(x)=\frac{(\alpha x)^2}{2}\quad x\in\mathbb{R},\; \alpha\in(0,\infty),
\end{equation}
in which the process reads 
\begin{equation}
    dZ(t)=BZ(t)\,dt + G\,dW(t),\quad Z(0)=(Q_0,p_0)^\top,
\end{equation}
where $B$ is defined such that $BZ(t)=(J-R)\nabla H(Z(t))$. 
In this case, the resulting process converges for $t\to\infty$ in distribution to a normal distribution with known expectation and covariance matrix.

\section{Results and Discussion}\label{sec5}
To introduce the port-Hamiltonian formalism into the neural network architecture, we can analyze the behavior of port-Hamiltonian Neural Networks (pHNN) on different tasks such as the mass-spring system with the damped term, external force and the Duffing equations. 
In particular, the port-Hamiltonian Neural Network (pHNN) is a significant advance in learning and predicting the dynamics of non-autonomous systems. 
Many real-world dynamical systems involve time-dependent forces and energy dissipation, which pose a challenge for learning. 
Desai et al.\ \cite{[10]:} evaluate the pHNN on several tasks, 
including a mass-spring system with damping and an external force, and a Duffing system. 
In particular, the pHNN can visually recover the Poincaré section of a chaotically driven system, 
highlighting its potential to identify and understand chaotic trajectories with minimal data. 
Its applicability extends to complex nonlinear forced and damped physical systems, 
encouraging applications in areas such as chemical bonding forces, robotic motion, and controlled dynamics without explicit knowledge of force and damping.

\begin{figure}[htbp]
\centering
\includegraphics[width=.90\textwidth]{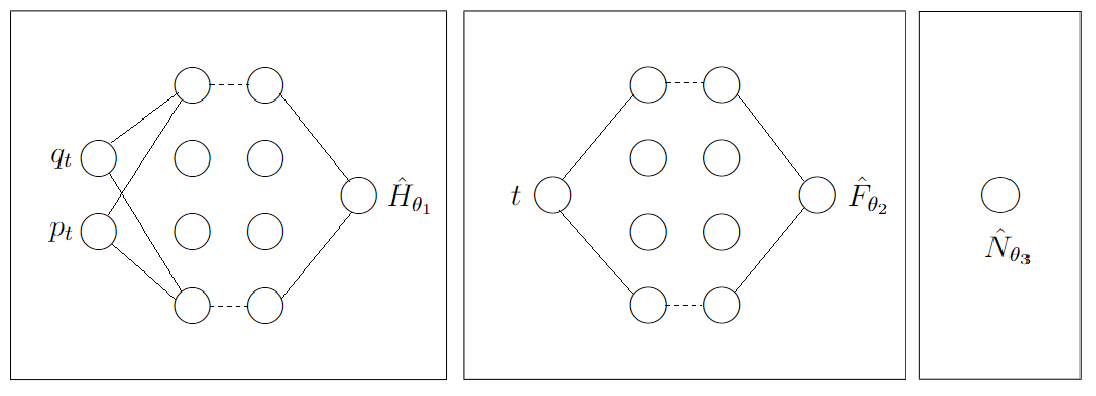}
\caption{Schematic representation of a port-Hamiltonian neural network (pHNN).}\label{figure4}
\end{figure}

As we can see in Figure~\ref{figure4}, the main idea is to use port Hamiltonian theory to explicitly learn the force $F_{\theta_2}$,
the damping term $N_{\theta_3}$ and the Hamiltonian $\mathcal{H}_{\theta_1}$. 
This should allow us to predict the time derivatives
\begin{equation}\label{mat}
    \begin{bmatrix}
        \hat{\dot{q}}_t\\
        \hat{\dot{p}}_t
    \end{bmatrix}=\begin{bmatrix}
        \frac{d\hat{H}_{\theta_1}}{dp_t}\\
        -\frac{d\hat{H}_{\theta_1}}{dq_t}+\hat{N}_{\theta_3}\frac{d\hat{H}_{\theta_1}}{dp_t}+\hat{F}_{\theta_2}
    \end{bmatrix}.
\end{equation}
Furthermore, knowing $[\dot{q},\dot{p}]$ from the data, we can compute the loss function \cite{Wang}
\begin{equation}\label{loss}
    \mathcal{L}=\underbrace{||\hat{\dot{q}}_t-\dot{q}_t||_2^2+||\hat{\dot{p}}_t-\dot{p}_t||_2^2}_{\text{first part}}
    +\underbrace{\lambda_F||\hat{F}_{\theta_2}||_1+\lambda_N||\hat{N}_{\theta_3}||_1}_{\text{second part}}.
\end{equation}
In particular, the procedure for implementing the port-Hamiltonian neural network is as follows:
\begin{itemize}
    \item Obtain ground truth state variable data $[q, p, t]$ and time derivatives $[\dot{q}, \dot{p} ]$ from trajectories of a given system;
    \item Provide state variable information to pHNN, which learns a Hamiltonian, force, and damping term to predict $[\hat{\dot{q}}, \hat{\dot{p}}]$;
    \item Optimize pHNN by minimizing the loss function of equation \eqref{loss};
    \item Train the system;
    \item Use the pHNN in a scientific integrator to evolve a set of random initial conditions in the test set;
    \item Estimate the goodness of the performance by computing the mean square error of the predicted state variables.
\end{itemize}

To understand the power of this method, we again compare different neural network models: 
the baseline neural network, the Hamiltonian neural network, the time-dependent Hamiltonian neural network, and the port-Hamiltonian neural network, against ground truth data \cite{[10]:}. \\
In a simple mass-spring system, characterized by mass $m$ and spring constant $k$, the time-independent Hamiltonian is given by
\begin{equation*}
    H = \frac{1}{2} kq^2 + \frac{p^2}{2m}\,.
\end{equation*}
After training the system and testing it over 25 random initial conditions, 
we find that both the Hamiltonian Neural Network (HNN) and the port-Hamiltonian Neural Network (pHNN) produce similar results in terms of state and energy mean squared error (MSE). 
In particular, pHNN successfully learns to predict force and damping over time, 
even when the ground truth observations are null. 
However, the prediction error remains small, approximately on the order of $10^{-5}$ and $10^{-8}$, respectively. 
Furthermore, comparing the state and energy MSE among the listed methods, 
it is evident that HNN achieves the lowest MSE, followed by pHNN, time-dependent Hamiltonian neural network (TDHNN), and the baseline neural network (bNN).\\
In the damped mass-spring system, we introduce a damping term that gradually decreases the initial energy of the system over time. 
Denoted by $\delta$, the damping coefficient, the system is described by the equation
\begin{equation*}
  \dot{q} = -\delta q\,.
\end{equation*} 
Furthermore, based on the definition $p = m\dot{q}$, we have:
\begin{equation*}
   \dot{q} = \frac{1}{m}p\,.
\end{equation*} 
By integration, we get $q^\top q - \delta q^\top \dot{q}$, which, when combined with a kinetic energy term, gives the Hamiltonian:  
\begin{equation*}
    H = \frac{p^\top p}{2m} + \frac{q^\top q}{2} - \delta q^\top \dot{q}\,.
\end{equation*}  
However, this formulation violates the previous equation, since 
$\frac{\partial H}{\partial p} = \dot{q} - \delta mq^{-1} \dot{q}\neq\dot{q}$. 
Consequently, including this term violates energy conservation, making it impossible to derive a scalar Hamiltonian for such a system.
 After training and testing this system, it is observed that both the Hamiltonian Neural Network (HNN) and the Time Dependent Hamiltonian Neural Network (TDHNN) fail to capture the dynamics of the system accurately. 
However, the baseline neural network (bNN) and the \textit{port-Hamiltonian neural network} (pHNN) show better performance, effectively recovering the system dynamics and minimizing the state and energy mean square error (MSE). 
Moreover, in this scenario, the pHNN converges to forcing and damping terms consistent with the ground truth.\\ 
In the forced mass-spring system, we introduce a forced time-dependent force that governs the undamped system. 
A Hamiltonian describing this forced mass-spring system has the form 
\begin{equation*}
   H = \frac{1}{2}kq^2 + \frac{p^2}{2m} - qF_0 \sin(\omega t) \,.
\end{equation*}
Here, $F_0$ is the force amplitude and $\omega$ is the frequency of the external force term.
After training and testing, both the Hamiltonian neural network (HNN) and the time-dependent Hamiltonian neural network (TDHNN) have difficulty learning the dynamics due to the lack of explicit time dependence in the Hamiltonian. 
Conversely, the port-Hamiltonian neural network (pHNN) exhibits the lowest state, and energy mean squared error (MSE), indicating better performance. 
In addition, the baseline neural network (bNN) has a relatively low MSE. 
However, when we consider a more complex force term, such as 
\begin{equation*}
  H = \frac{1}{2}kq^2 + \frac{p^2}{2m} - qF_0 \sin(\omega t) \sin(2\omega t), 
  \end{equation*}
   where the complexity of the force significantly affects the performance of bNN.
   In this scenario, the state and energy MSE of bNN becomes comparable to that of HNN. 
   Meanwhile, pHNN continues to capture and evolve the non-harmonic force accurately, demonstrating superior performance in initial state tests compared to other models.

The problem with the presented neural networks is that in the real world, such neurons are noisy, and the output is a probabilistic input function. 
Therefore, in the appendix, we analyze some stochastic representations of the NNs. 

According to \cite{[53]:}, we consider the structure of a port-Hamiltonian system of the form
\begin{equation}
    \begin{cases}
        \dot{x}=&[J(x)-R(x)]\frac{\partial H}{\partial x}(x)+G(x)u,\\
        y=&G^\top(x)\frac{\partial H}{\partial x}(x),
    \end{cases}
\end{equation}
with $x\in\mathbb{R}^n$ the energy variable, $H(x):\mathbb{R}^n\to\mathbb{R}$ the total stored energy, $u,y\in\mathbb{R}^m$ the port power variables.\\
In this context, we can explore Evolution Strategies (ES) as a stochastic optimization technique to evolve a population of solutions. 
ES includes basic concepts such as mutation, selection, and fitness evaluation. These strategies simulate collective learning processes within a population of individuals.

ES is part of a broader category of evolutionary algorithms, including genetic algorithms (GA) and evolutionary programming. 
These methods start with a random population and evolve through selection, mutation, and recombination processes. 
The interaction between creating new genetic information and its evaluation and selection drives evolution. 
Individuals who perform better have an increased chance of survival and produce offspring that inherit their genetic information.

In ES, the control synthesis problem is simplified to finding the minimum of a function, called fitness, within a feasible set of values.
 Each individual in the population is characterized by a set of exogenous parameters $\pi$ and endogenous parameters $\sigma$ representing points in the search space and independent paths in the phase space, respectively. 
Sufficient initial conditions are provided to avoid overfitting and excessive computational cost.

The goal is for the set of trajectories to approach an equilibrium point under the influence of a stabilizing controller. 
The fitness metric evaluates the performance of the controller, which is derived from the energy function $H_a$. 
Mutation and self-adaptation act independently on each individual, generating new parameters accordingly:
\begin{equation*}
  \pi_{i+1} = \pi_i + \sigma_i \cdot N(0, 1) 
\end{equation*}
\begin{equation*}
\sigma_{i+1} = \sigma_i \cdot \exp(\tau \cdot N(0, 1)) + \tau' \cdot N(0, 1)\,.
\end{equation*}
Here, $\tau$ and $\tau'$ are proportional to $(\sqrt{2n_\pi})^{-1}$ and 
$(\sqrt{2\sqrt{n_\pi}})^{-1}$, respectively, where $n_\pi$ represents the dimension of the search space. These equations govern the evolution of the parameters within the ES framework and facilitate the exploration of the solution space.

\section{Conclusions}\label{sec6}
The authors of \cite{[53]:} introduced a method that uses a neural network as an approximator to generate a gradient field for controller design, thereby avoiding the need to solve partial differential equations. 
To extend these results to a stochastic port-Hamiltonian scenario, it's necessary to introduce stochastic elements into the existing framework. 
This can be done by the following steps:
\begin{enumerate}
    \item Define a stochastic version of the port-Hamiltonian system by introducing stochastic differential equations (SDEs) to model the system dynamics. For example, consider a system described by
\begin{equation*}
dx = \bigl(J(x) - R(x)\bigr)\, \nabla H(x)\, dt 
    + G(x) u \,dt + \Sigma(x)\, dW, 
\end{equation*}
Where $\Sigma(x)$ is a matrix characterising the intensity of stochastic perturbations and $dW$ is a standard Wiener process.
    \item Extend the control synthesis methodology to account for the stochastic nature of the system. Modify the control law to ensure stability and performance despite stochastic disturbances. The control law might take the form
   \begin{equation*}
  u = -K(x) + \mu(x)
\end{equation*}
where $K(x)$ is derived from the deterministic part of the system (as in the original paper) and $\mu(x)$ is an additional term to counteract stochastic disturbances.
\end{enumerate}

Accordingly, an appropriate Lyapunov function must be identified or constructed to ensure the stability of the stochastic system. Lyapunov analysis for stochastic systems typically uses It\^o's formula and aims to show that the expected value of the Lyapunov function decreases with time. Using evolutionary strategies, adapt the optimization process to handle the stochastic nature of the system. This may involve modifying the fitness function to account for stochastic behavior, possibly by considering expected performance over a range of stochastic disturbances.
 Finally, a new case study will extend the ball and beam system to include stochastic elements. Simulations can be used to demonstrate the stochastic disturbances.

\appendix

\section{Stochastic Neural Networks}\label{appA}
Neurons, the basic computational units in the brain, are known to form neural networks (NNs) by connecting through synaptic junctions. 
However, real neurons introduce noise that makes the output a probabilistic function of the input. Unlike deterministic NNs, the activation map is stochastic. 
These networks are better at avoiding local minima during training. 
They are better suited to tasks with noisy or incomplete data, although implementation and training complexity increases in the stochastic case. 
In particular, their randomness improves their robustness against overfitting since they do not learn the noise in the training data precisely as deterministic networks do. 

\begin{figure}[htbp]
\centering
\setlength{\unitlength}{0.8cm}
\includegraphics[width=0.7\textwidth]{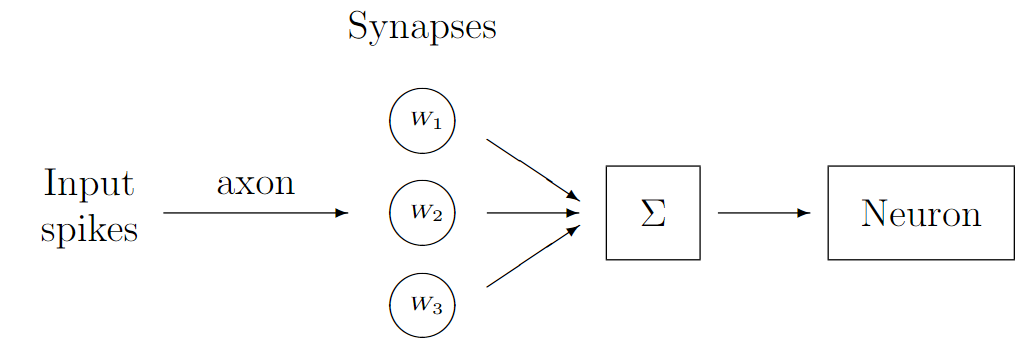}
\caption{Schematic representation of a basic artificial Neuron 
block.}\label{figure18}
\end{figure}

The interconnection of neurons from different layers creates a neural network. 
The primary components for the synapses and neurons that make up a stochastic neural network are the stochastic neurons implemented by magnetic tunnel junctions (MTJs). 
Without going into the physical details, it can be described that an MTJ device consists of two ferromagnets separated by a thin insulator and is characterized by a switching (activation) probability (see Zhu and Park~\cite{[44]:} for more technical details). 
In the artificial neural network model, input spikes flow into synapses, each assigned a weight, as shown in Figure~\ref{figure18}.

Then  
\begin{equation*}
    \Sigma:=\text{\textit{Input}}\times\text{\textit{Weights}}=\text{\textit{Summed Output}}
\end{equation*}
defines the weighted summation and finally the neuron fires depending on a threshold or activation function. 
Vreeken~\cite{[45]:} introduced \textit{Spiking Neural Networks} (SNNs), artificial neural networks designed to mimic natural neural networks closely. 
Neurons emit these short bursts of electrical energy when they accumulate in sufficient numbers. 
The spike travels along the neuron's axon and is controlled by the synapse, which consists of the end of the axon, a synaptic gap, and the initial part of the dendrite. 
Neurons must then use the spatial and temporal information of incoming spike patterns to encode their message to other neurons. 

A modification of stochastic neural networks is presented by Yu et al.~\cite{[47]:} with a \textit{Simple and Effective Stochastic Neural Network} (SE-SNN). 
This approach is based on the notion that activation uncertainty can be modeled at each layer by predicting a Gaussian mean and variance and then sampling the layer during the forward pass. 
This method emphasizes the activation distribution and incorporates an activation regularizer to optimize models with high uncertainty while still predicting the target label. 
SE-SNN exhibits favorable properties in pruning, adversarial defense, learning with label noise, and improving model calibration.
 In the context described, instead of encoding data using binary numbers, we can consider the representation by the probability of encountering 1's in bit streams. 
 Stochastic Computing (SC) uses this technique to represent continuous values using streams of random bits. 
 It uses conventional digital logic to perform computations based on stochastic bit streams.
SC, as described in \cite{[42]:}, encompasses a wide range of techniques, but the focus here is on a specific aspect used in \cite{[31]:}. 
The main benefits of SC include reducing hardware complexity and enabling fault-tolerant computing. This reduction in complexity comes from implementing functions such as the sigmoid, hyperbolic tangent, and exponential functions using linear finite state machines, resulting in lower hardware costs. 
However, this may result in reduced computational accuracy.
 In addition, as noted in \cite{[31]:}, the noise introduced in SC can mitigate overfitting problems and improve inference accuracy. 
 Consequently, both input and output are represented by bit streams, with values encoded as probabilities of encountering 1's in these streams. 
 This approach emphasizes the fusion of probabilistic methods with digital computation, providing a unique data representation and processing paradigm.



\end{document}